\documentclass{amsart}

\usepackage{amsmath}
\usepackage{amssymb}

\def \bbF {{\mathbb F}} 
\def \bbZ {{\mathbb Z}}
\def \bbQ {{\mathbb Q}}

\def \Ord {\mathcal{O}} 
\def \calp {\mathfrak{p}}

\begin{document}

\title{Computing quadratic subfields of number fields}

\author{Andreas-Stephan Elsenhans and J\"urgen Kl\"uners}

\begin{abstract}
Given a number field, it is an important question in algorithmic number theory to determine all its subfields.
If the search is restricted to abelian subfields, one can try to determine them by using class field theory.
For this, it is necessary to know the ramified primes.
We show that the ramified primes of the subfield can be computed efficiently. Using this information
we give algorithms to determine all the quadratic and the cyclic cubic subfields of the initial field. 
The approach generalises to cyclic subfields of prime degree.
In the case of quadratic subfields, our approach is much faster than other methods.
\end{abstract}

\maketitle

\section{Introduction}

In recent years several people have worked on improved methods to
compute the subfields of a given number field
$L$~\cite{EK,Kl,KHN,HS1}. Depending on the situation the more
combinatorial approach used in \cite{Kl} or the one based on
generating subfields introduced in \cite{KHN} might be preferable. The
differences of these approaches are analised in \cite{EK}. In the latter article
we also introduce the notation composition-generating. Here we partially answer the
question if it is possible to efficiently compute small degree subfields. We remark
that in the hardest examples with elementary abelian Galois group of 2-power degree
we compute such a set within a few seconds for a degree 128 extension.

In this note we present an approach for quadratic or more generally cyclic subfields. 
The main idea is that the ramified primes of a subfield are ramified in $L$. As soon 
as the ramified primes are known only finitely many fields are left. 
In principle the ramified primes can be determined by factoring the discriminant of 
the defining polynomial. But this is usually very costly as the discriminant can be very large. 

The main new idea is to compute a by far smaller number that contains all 
the odd ramified primes of each quadratic subfield.

We implemented the new approach in {\tt magma}. Further, we compare the run time with other implementations. 

\subsection*{Notation}
For a monic polynomial $f$ we call a monic polynomial $g$ such that $f = g^e$ an $e$-th root of $f$.

\section{Fields, subfields and ramification}

\subsection{Lemma}
Let the field $L = \bbQ[X] / (f(X))$ be given by a monic polynomial
$f \in \bbZ[X]$.  Assuming the field $K$ of degree $e$ is a proper
subfield of $L$ then there is a monic polynomial $g \in \Ord_K[X]$
such that $f = N_{K/\bbQ}(g)$.

Further, let $(p) = \calp^e$ be a totally ramified prime with a prime ideal $\calp \subset \Ord_K$. Then the reduction of
$f$ modulo $p$ is an $e$-th power.

{\bf Proof:} As $p$ is totally ramified the reduction of $N_{K/\bbQ}(g)$ modulo $\calp$ is an $e$-th power. This reduction coincides with the
reduction modulo $p$. \qed

\subsection[Power detection]{Remark}
\begin{enumerate}
\item
Let a monic polynomial  $f \in R[X]$ of degree $e \cdot n$ be given. 
We have to check if $f$ is an $e$-th root. For this we assume $f = (X^{en} + a_{en -1}X^{ne-1} + \cdots + a_1 X + a_0)
= g^e = (X^n + b_{n-1} X^{n-1} + \cdot + b_1 X + b_0)^e$.
A comparison of coefficients shows that the 
coefficients $b_i$ are uniquely determined by $a_{en-1},\ldots,a_{(e-1)n}$ as long as $e$ is a unit in $R$.
\item
In case  $f \in \bbQ[X]$ Newton's identities can be used to compute the power sums of the roots 
$r_1,\ldots,r_{en}$ of $f$ without using the roots. I.e.
\begin{eqnarray*}
s_1 &:=& r_1 + \cdots + r_{en} \\
s_2 &:=& r_1^2 + \cdots + r_{en}^2 \\
&\vdots&
\end{eqnarray*} 
Using Newton's identities a second time to compute the polynomial to the power sums of roots 
$\frac{s_1}{e},\frac{s_2}{e},\ldots,\frac{s_n}{e}$ results in the polynomial $g$ as well.

We refer to~\cite{BFSS} for efficient implementations of Newton's identities and further use of them.
\end{enumerate}

\subsection{Algorithmic approach}
A number field  $L = \bbQ[X] / (f(X))$ with a cyclic subfield $K$ of degree $e$. 
Searching for tamely ramified primes in $K$ can be done as follows:
\begin{enumerate}
\item 
Computing a candidate $g$ for the $e-$th root of the polynomial $f$. By using only the terms of degree $X^{en},\ldots,X^{(e-1)n}$ of $f$.
\item
Compute the gcd of the numerators of all the coefficient of $f - g^e$.
\item
The prime divisors of the gcd are the only possibly purely tamely ramified primes of $K$.
\end{enumerate}

When the tamely ramified primes are found, we add the possibly wildly ramified primes, i.e. the primes dividing $e$. 
If $K$ can be an imaginary field, we add the prime $\infty$ respectively $-1$ as a possible divisor of the discriminant of $K$.

In practice the number that is factored in step 3 is usually small. The total factorization time for all our test polynomial is about 
0.03 seconds.

\subsection{Determination of all quadratic subfields}
Let a number field $L$ of even degree be given. Use the approach above to compute a set of potential factors (prime factors and $-1$)
of the discriminant $\Delta$ of a quadratic subfield $\bbQ(\sqrt{\Delta}) \subset L$. 

We can enumerate the finitely many possible values of $\Delta$ that are compatible with this by multiplying the elements of each subset of 
the set of potential divisors of $\Delta$. For each candidate $\Delta$ we can test directly if $X^2 - \Delta$ has a root in $L$.
This gives all the quadratic subfields.

\subsection{Improvement by cycle types}
Given a number field $L = \bbQ[X]/(f(X))$ by a monic polynomial
$f \in \bbZ[X]$ and a prime number $p$ that is unramified in every
quadratic subfield.  If the $p$-adic factorization of $f$ has a factor
of odd degree that $p$ is split in each quadratic subfield of $L$.  If
it is impossible to select some of the factors of this $p$-adic
factorization in such a way that their product has half the degree of
$f$ then $p$ is inert in each quadratic subfield of $L$. If a prime is
inert in two subfields of degree 2, it must be split in a third
one. Therefore this case means that there is at most one subfield of
degree 2.

If a prime is forced to be at the same time inert and split in all quadratic subfields, it is shown that no quadratic subfield exists. 

If we write the discriminant $\Delta = (-1)^{e_0} p_1^{e_1} \cdots p_k^{e_k}$ with $e_0,\ldots,e_k \in \{0,1\}$ and 
potentially ramified finite primes 
$p_1,\ldots,p_k$ each prime with known decomposition in $\bbQ(\sqrt{\Delta})$ gives a $\bbF_2$-linear equation for $e_0,\ldots,e_k$. 
By inspecting the $p$-adic factorization of $f$ for several primes $p$ we can generate a $\bbF_2$-linear system for $e_0,\ldots,e_k$. 
This reduces the number of quadratic fields that have to be tested.

In the vast majority of test examples, all remaining subfield candidates are in fact subfields. 

\subsection{Improvement by twists of subfields}
As soon as the subfield $\bbQ(\sqrt{\Delta_1})$ and $\bbQ(\sqrt{\Delta_2})$ are found, the subfield $\bbQ(\sqrt{\Delta_1 \Delta_2})$
is for free. Further, if the subfield $\bbQ(\sqrt{\Delta_1})$ is known and $\bbQ(\sqrt{\Delta_2})$ is excluded, 
$\bbQ(\sqrt{\Delta_1 \Delta_2})$ is excluded as well. One can easily use these facts to further reduce the number of root computations
of $X^2 - \Delta$  in $L$.

For example, let $L / \bbQ$ be a degree 128 extension with Galois group $C_2^7$. Then $L$ has 127 quadratic subfields. We need do determine
7 quadratic subfields directly. All the others can be constructed out of them by twisting. 

\section{Generalization - cyclic cubic fields}

\subsection{Subfields as Kummer extensions}
With the methods described above we can determine a set containing all ramified primes 
of all cyclic cubic subfields of a given field. By Kummer theory a cyclic cubic extension 
can be described easily as radical extension of $\bbQ(\zeta_3)$.
Thus, the field 
$$
K_0 := \bbQ(\zeta_3)(\sqrt[3]{\zeta_3},\sqrt[3]{\pi_1},\ldots,\sqrt[3]{\pi_k})
$$
with $\pi_1,\ldots,\pi_k$ the primes of $\bbZ[\zeta_3]$ above the potentially ramified primes 
contains all cyclic cubic extensions of the given field.

As above we can determine primes that split in all cyclic cubic subfields. By using 
$\bbF_3$ linear algebra we can replace $K_0$ by a subfield $K_1$. Next we have to list all 
the cyclic cubic subfields of $K_1$.

\subsection{Cubic extensions of $\bbQ$ from Kummer extensions}
Assume that the extension $K_2 := \bbQ(\zeta_3,\sqrt[3]{a}) / \bbQ$ with $a \in \bbQ(\zeta_3)$ 
is abelian of degree 6.
Then the norm $N(a)$ has to be a cube in $\bbQ$. 
To compute  the cubic subfield, we denote 
by $\sigma \colon \bbQ(\zeta_3) \rightarrow \bbQ(\zeta_3)$ the non-trivial automorphism. 
Then the image of $\sqrt[3]{a}$ under the automorphism of order 2 of $K_2$ is given by
$$
g := \left(\frac{\sqrt[3]{a}}{\sqrt[3]{N(a)}}\right)^2  \sigma(a) .
$$
Thus, a symmetric expression such as $g + \sqrt[3]{a}$  will generate the cubic subfield. We remark that
a more general approach is described in \cite{Fi}.
After a determination of the minimal polynomial of the generator found, we can compute its roots in the initial
field $L$ to check that $\bbQ(g)$ is a subfield of $L$.

\subsection{Remark}
In principle it would be possible, to determine cyclic subfields of any degree in this way. 
But, as the degree of the cyclotomic fields that have to be used get larger, the approach
will be less efficient.

\section{Performance tests}
We implemented the above approach in {\tt magma}\cite{BCP} and compared the performance with {\tt magma's} subfield algorithm 
restricted to a search for quadratic subfields.
As a test sample we used on polynomial for each transitive group in degree 16, 18 and 20. (In total 1954 + 983 + 1117 = 4054 polynomials.) 
Computing the quadratic subfields with our approach took 15 / 4 / 15 seconds. Using {\tt magma} we can compute them in 24 / 17 / 35 seconds. 

Further, some test polynomials are listed in~\cite{HS2}.  
The largest example is a degree 128 field with Galois group $(\bbZ /2 \bbZ)^7$ and  127 quadratic subfields.
Computing generators of the lattice of all subfields takes about 50 minutes (48 minutes for LLL-reduction) in {\tt magma}. 
The time to compute all the subfields is listed in~\cite{HS1} as 86 minutes.
We can determine the quadratic subfields in 22 seconds. In this example the quadratic subfields generate the lattice of all 
subfields, i.e. this set is composition generating.

Further, we can determine the quadratic subfields of all other 23 examples in the list~\cite{HS2} in 3.6 seconds. 
Only two of them are listed in~\cite{HS1} with less than 20 seconds of CPU time.

\medskip

Further, we tested the computation of cyclic cubic subfields in the same way. 
We took the same 983 polynomials of degree 18. Computing all cubic subfields with {\tt magma} took 22 seconds.
Using our approach, we find all the cyclic ones in 5.9 seconds. When we restrict to those polynomials  
with at least one such subfield, {\tt magma} takes 6.2 seconds whereas the above method takes 4.2 seconds.

Within the examples of~\cite{HS2}, the hardest one is a degree 81 polynomial with a Galois group of 
order 162 that results in 4 cubic subfields. We can determine them in 1.9 seconds. 
Computing all 35 subfields in {\tt magma} takes 117 seconds.
Here, {\tt magma} uses 109 seconds for LLL-reduction
to find generators of the subfield lattice.
In~\cite{HS1} the run time for all subfields of this example is listed as 716 seconds. 

\medskip

All computations are done by using {\tt magma 2.24-7} on one core of an Intel i7-7700 CPU running at 3.6\,GHz. 
The code is available on the authors web pages.

\end{document}